\documentclass[12pt]{amsart}
\usepackage{bbm}
\usepackage{amsmath} 
\usepackage[english]{babel}
\usepackage{color}
\usepackage{anysize}
\usepackage{geometry}
\marginsize{2.0cm}{2.0cm}{2.0cm}{2.0cm}
\usepackage{url}
\usepackage{hyperref}
\usepackage[normalem]{ulem}
\usepackage{mathtools}
\mathtoolsset{showonlyrefs}
\usepackage{soul}
\usepackage{url}
\mathtoolsset{showonlyrefs}
\usepackage{graphicx}
\usepackage{tikz}
\usepackage{atbegshi}
\theoremstyle{plain}
\numberwithin{equation}{section}
\newtheorem{theorem}{Theorem}[section]
\newtheorem{corollary}[theorem]{Corollary}
\newtheorem{proposition}[theorem]{Proposition}
\newtheorem{lemma}[theorem]{Lemma}
\newtheorem{remark}[theorem]{Remark}
\newcommand{\rmd}{{\rm d}}
\newcommand{\qedw}{\hfill \ensuremath{\Box}}
\title[Wasserstein error between telegraph process and Brownian motion]{Wasserstein error estimates between telegraph processes and Brownian motion}
\author{Gerardo Barrera}
\address{Center for Mathematical Analysis, Geometry and Dynamical Systems, Mathematics Department, Instituto Superior T\'ecnico, Universidade de Lisboa, 1049-001, Av. Rovisco Pais 1, 1049-001 Lisboa, Portugal.
\url{https://orcid.org/0000-0002-8012-2600}}
\email{gerardo.barrera.vargas@tecnico.ulisboa.pt}

\author{Jani Lukkarinen}
\address{
Department of Mathematics and Statistics,
University of Helsinki.
P.O. Box 68, Pietari Kalmin katu 5, FI-00014, 
Helsinki, Finland.
\url{https://orcid.org/0000-0002-8757-1134}}
\email{jani.lukkarinen@helsinki.fi}
\author{Mikko S. Pakkanen}
\address{
Department of Mathematics, Imperial College London.
South Kensington Campus
London SW7 2AZ, London, United Kingdom.
\url{https://orcid.org/0000-0002-0696-4914}}
\email{m.pakkanen@imperial.ac.uk}

\begin{document}

\AtBeginShipoutFirst
{
\begin{tikzpicture}[remember picture, overlay]
\node[anchor=north west, xshift=1.5cm, yshift=-2.5cm] at (current page.north west){
\includegraphics[width=4cm]{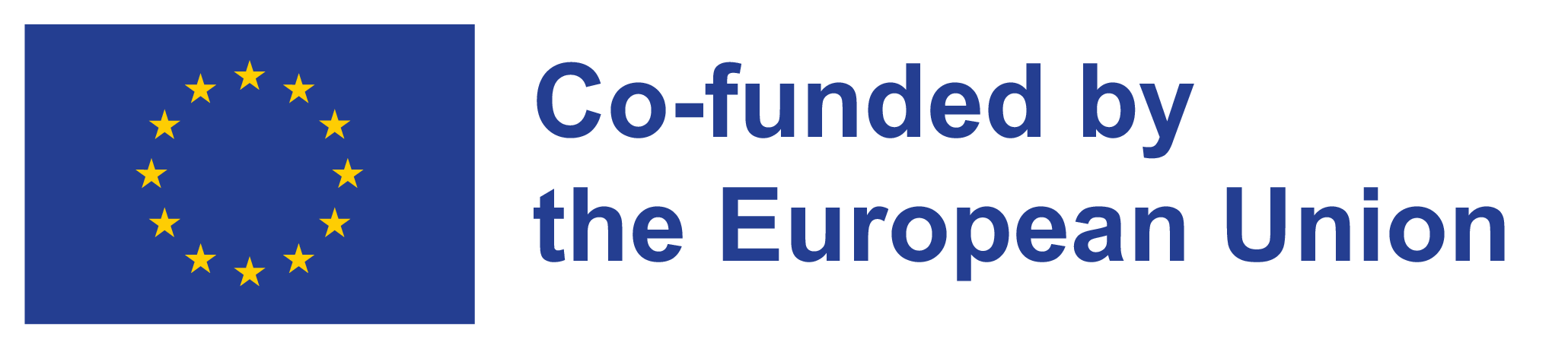}  
};
\end{tikzpicture}
}

\begin{abstract}
We provide non-asymptotic error bounds in the path Wasserstein distance with quadratic integral cost between suitable functionals of the telegraph process and the corresponding  functional of Brownian motion with explicit diffusivity constant.
These results cover, in particular, the well-known example of the exponential integral functional of the geometric Brownian motion. 
The non-asymptotic error bounds tend to zero in the so-called Kac regime.
Moreover, the error bounds
remain valid when the flip rate for the telegraph process is small.
We assess the sharpness of the error bounds through numerical experiments.
\end{abstract}

\maketitle


\section{\textbf{Introduction, the model, main results and consequences}}\label{sec:intr}

\subsection{\textbf{The model}}
\hfill

\noindent
The so-called Goldstein--Kac (symmetric) telegraph process  is one of the simplest examples of a random evolution, 
see~\cite{Goldstein1951} 
and~\cite{Kac1974}.
It is a mathematical description of non-interacting particles moving in one dimension with alternating (back and forth) two  finite velocities.
It has been proposed as an alternative model for diffusion models such as the celebrated Black--Scholes--Merton model based on geometric Brownian motion.
Therefore, it has been extensively considered and studied by different communities in applied science such as the probability, finance, biology, ecology and physics communities. 
It is an important and well-studied mathematical object in its own right~\cite{Bogachev2011,Cinque2023,DeGregorio2021,
DiCrescenzo2013,DiCrescenzo2018,
Iacus2009,
Janssen1990,
Kolesnik2018,Kolesnik2015,
KolesnikTurbin1998,
LopezRatanov2014spl,LopezRatanov2014,
Martinucci2020,Ratanov2020-2,Ratanov2021-2}, and it has broad generalizations with vast applications to physical systems,
ecological models, biological 
models~\cite{DeGregorio2021-2,DiCrescenzo2023oct,Foong1994,
Ghosh2014,Goldstein1951,Kac1974,
Kolesnik1998,Orsingher1990,Ratanov2019,
Stadje2004}, finance and risk models~\cite{DiMasi1994,Kolesnik2013,LopezRatanov2012,Mazza2004,
Pogorui2021,Ratanov2007,Ratanov2009,Ratanov2016}, etc. Further references may be found from therein.

In the sequel, we introduce the Goldstein--Kac (symmetric and asymmetric) telegraph processes.
More precisely, it describes the movement of a particle which starts at time zero
from the origin and moves with a finite non-zero constant velocity $v_0$  on the line with the following rule:
\begin{itemize}
\item[(i)] the initial direction of the motion ($v_0$ or $-v_0$) is chosen at random with the uniform probability $\mu$, that is, $\mu(\{v_0\})=\mu(\{-v_0\})=1/2$, 
\item[(ii)] the dynamics of its direction is
driven by a homogeneous Poisson clock (process) of a positive constant intensity $\lambda$. In other words, 
when the Poisson clock rings, the particle instantaneously takes the opposite direction and keeps moving with the same velocity until the next ring in the Poisson clock happens, then it takes the opposite direction again, and so on. Recall that
the inter-arrival random times for a homogeneous Poisson process of intensity $\lambda$ are independent and identically distributed random variables having Exponential law with parameter $\lambda$.
\end{itemize}
The latter defines a random process $X:=(X(t))_{t\geq 0}$ in a probability space $(\Omega,\mathcal{F},\mathbb{P})$.
We note that $|X(t)|\leq |v_0|t$ for all $t\geq 0$.
The absolute continuous part of particle position $p(x,t)$ spread over the interval $x\in (-|v_0|t, |v_0|t)$, formally defined as $p(x,t)=\mathbb{P}(X(t)\in \rmd x)/\rmd x$,
solves the  hyperbolic second-order linear partial differential equation (a.k.a. damped wave equation)
\begin{equation}\label{eq:pde1}
\frac{\partial^2}{\partial{{t}^2}} p(x,t)+2\lambda \frac{\partial}{\partial_{{t}}} p(x,t)=v^2_0 \frac{\partial^2}{\partial{{x}^2}} p(x,t)\quad \textrm{ with }\quad p(x,0)=\delta_{0}(x),
\end{equation}
where $\delta_0(x)$ denotes the Dirac delta function with total mass at zero. 
We point out that the solution of~\eqref{eq:pde1}, a.k.a. Green's function, can be represented in terms of Bessel functions. 
It may be used as a model of the one-dimensional transport produced by diffusion with finite propagation speed.
We refer to~\cite{Foong1994},~\cite{Kabanov1992},~\cite{Kac1974},~\cite{Orsingher1990} and Chapter~3 in~\cite{Kolesnik2013} for more details.

In the sequel, we introduce the generalized telegraph process.
Let $N:=(N(t))_{t\geq 0}$ and $N^*:=(N^*(t))_{t\geq 0}$ be two independent Poisson processes with positive intensities (rates) $\lambda$ and $\lambda^*$, respectively.
Let $\mathcal{E}:=\{v_0,-v^*_0\}\subset \mathbb{R}$ be a set of two elements (space state),
 and equip it with the uniform probability measure $\mu$, that is, $\mu(\{v_0\})=\mu(\{-v^*_0\})=1/2$.
We consider the \textit{velocity process} $V:=(V(t))_{t\geq 0}$  taking values in the state space $\mathcal{E}$ with C\`adl\`ag (continue \`a droite, limite \`a gauche) paths $t\mapsto V(t)$   described as follows: 
\begin{itemize}
\item[(1)] choose at random $v_0$ or $v^*_0$ according to $\mu$ and independent of the Poisson processes $N$ and $N^*$,
\item[(2)] then remain in this state according to the Poisson process with the same label, in other words, if we choose $v_0$ in Item~(1) then $V(t)=v_0$ for all times $t$ up to the corresponding Poisson clock $N$ rings, otherwise,  $V(t)=v^*_0$ for all times $t$ up to the corresponding Poisson clock $N^*$ rings,
\item[(3)] when the Poisson clock rings, change to the state with the other label and remain in this state according to its Poisson clock,
\item[(4)] the switching between the states is repeated and therefore the velocity process $V$ is well-defined for all no-negative times due to the well-known properties of the Poisson processes.
\end{itemize}
Then we define the (generalized) \textit{telegraph process} $X:=(X(t))_{t\geq 0}$ as the position process starting at $z\in \mathbb{R}$ of the velocity process, that is,
\begin{equation}
X(t):=z+\int_{0}^t \rmd s\, V(s)\quad \textrm{ for any }\quad t\geq 0.
\end{equation}
We point out that $z-\max\{|v_0|,|v^*_0|\}t\leq X(t)\leq z+\max\{|v_0|,|v^*_0|\}t$ for all $t\geq 0$, where $\max\{x,y\}$ denotes the maximum between the real numbers $x$ and $y$.
The absolutely continuous component of the particle position solves a hyperbolic second-order linear partial differential equation
whose explicit solution can be written in terms of Bessel functions, see 
Section~2~in~\cite{Beghin2001}.

For shorthand notation we always assume that $z=0$.
We point out that the multicomponent process $(X,V)$ is Markovian with respect to the natural filtration, however, the projection $X$ is no longer Markovian.

The different velocities of motion, i.e., $v_0+v^*_0\neq 0$,  naturally brings asymmetry in the sample paths $t\mapsto X(t)$ and hence in many observables of them. This makes the evaluation of the distribution for observables \textit{a priori}
much more difficult than in the classical symmetric case.
When $\lambda=\lambda^*>0$ and $v_0=-v^*_0$ we say that the telegraph process is symmetric and we denote it by $X^{\textsf{sym},v_0}=(X^{\textsf{sym},v_0}(t))_{t\geq 0}$.

\subsection{\textbf{The weighted Lusin--Lipschitz observables and the path Wasserstein distance}}
\hfill

\noindent
In this subsection, we 
define the path space and the Lusin--Lipschitz observables that we are interested in, and
review the necessary background about Kantorovich--Rubinstein--Wasserstein path distance.

Let $T>0$ be fixed and set
\begin{equation}\label{eq:L2}
L^2_{\textsf{Ave}}([0,T],\mathbb{C}):=\left\{X:[0,T]\to \mathbb{C}: \frac{1}{T}\int_{0}^T \rmd s \,|X(s)|^2<\infty\right\},
\end{equation}
where $|\cdot|$ denotes the complex modulus.
For short, we write $L^2_{\textsf{Ave}}$ in place of $L^2_{\textsf{Ave}}([0,T],\mathbb{C})$.
We stress that when it is needed we write $L^2_{\textsf{Ave},\mathbb{R}}$ in place of $L^2_{\textsf{Ave}}([0,T],\mathbb{R})$.
Now, consider the set of continuous functions $C(L^2_{\textsf{Ave}},I)$ defined in $L^2_{\textsf{Ave}}$ and taking values in $I$, where $I$ may be chosen to be the complex numbers or the non-negative real numbers.
Let  $F\in C(L^2_{\textsf{Ave}},\mathbb{C})$ be an observable 
for which there exist  a non-negative constant $\kappa$ and  a weight function $W\in C(L^2_{\textsf{Ave}},[0,\infty))$ satisfying
\begin{equation}\label{eq:wlip}
|F(X)-F(Y)|\leq \kappa \left(\frac{W(X)+W(Y)}{2}\right) \|X-Y\|_{\textsf{Ave}}\quad \textrm{ for any }\quad X,Y\in L^2_{\textsf{Ave}},
\end{equation}
where
\[
\|X-Y\|_{\textsf{Ave}}:=\left(\frac{1}{T}\int_{0}^T \rmd s\, |X(s)-Y(s)|^2\right)^{1/2}
\]
is an average cost function.
We say that $F$ satisfying~\eqref{eq:wlip} is a ``$(\kappa,W)$-Lusin--Lipschitz" observable. We point out that the Lusin--Lipschitz observable~\eqref{eq:wlip} mimics the Lusin--Lipschitz inequality given in Remark~2.3 of~\cite{Brezis2020}.
For a continuous function $X:[0,T]\to \mathbb{R}$, for Lipschitz functions 
$f:\mathbb{R}\to \mathbb{R}$  and $g:\mathbb{R}\to \mathbb{R}$, and real constants $a$ and $b$, 
the following observable 
\begin{equation}\label{ec:nuev}
F(X):=g\left(\frac{1}{T}\int_{0}^T \rmd s \,f\big(e^{aX(s)+bs}\big)\right),
\end{equation} 
is Lusin--Lipschitz,
see Lemma~\ref{lemF} below for the precise form of $\kappa$ and $W$. 
In particular, the path $X$ can be taken to be random as a Brownian path, $f(z):=\ell z$, $z\in \mathbb{R}$ and  $g(z):=\max\{\ell z,c\}$, $z\in \mathbb{R}$ for some $\ell\in \mathbb{R}$. 
For $f$ and $g$ being the identity,~\eqref{ec:nuev}~becomes an integral exponential functional of the geometric Brownian motion, which has been widely used in financial modeling. It is known that the distribution of such integral exponential functional is non-trivial and exhibits an oscillating behavior which is related with the so-called Hartman--Watson distribution, see~\cite{Barrieu2004,HartmanWatson} and~\cite{Nandori2022} for a numerical evaluation.
We emphasize that an explicit representation for the law
of~\eqref{ec:nuev} is typically a difficult task, see~\cite{Kac1949} or Subsection~1.2 in~\cite{Sznitman1998}.

Let equip the product space $L^2_{\textsf{Ave}}\times L^2_{\textsf{Ave}}$ with its Borel product $\sigma$-algebra $\mathcal{B}(L^2_{\textsf{Ave}}\times L^2_{\textsf{Ave}})$ and 
consider the set of probability measures $\mathcal{P}_T$ in the measurable space 
\[\mathcal{J}_T:=(L^2_{\textsf{Ave}}\times L^2_{\textsf{Ave}},\mathcal{B}(L^2_{\textsf{Ave}}\times L^2_{\textsf{Ave}})).\]
For any $\mu,\nu\in \mathcal{P}_T$ we say that a probability measure ${\Pi}_*$ defined in $\mathcal{J}_T$ is a {\em coupling} between $\mu$ and $\nu$ if the marginals of ${\Pi}_*$ are $\mu$ and $\nu$, respectively. To be more precise,  
for any $B\in \mathcal{B}(L^2_{\textsf{Ave}})$ it follows that
\[
{\Pi}_*(B \times L^2_{\textsf{Ave}}) = \mu(B)\quad \textrm{ and }\quad
{\Pi}_*(L^2_{\textsf{Ave}} \times B) = \nu(B).\]
Let $\mathcal{C}(\mu,\nu)$ be the set of all couplings between $\mu$ and $\nu$.  
For any $\mu,\nu\in \mathcal{P}_T$, the average Wasserstein distance of order $2$ between $\mu$ and $\nu$, $\mathcal{W}_2(\mu,\nu)$, is defined by 
\begin{equation}\label{def:wa}
\begin{split}
\mathcal{W}_2(\mu,\nu):&=\inf\limits_{\Pi_* \in \mathcal{C}(\mu,\nu)}\left(\int_{L^2_{\textsf{Ave}} \times L^2_{\textsf{Ave}}} {\Pi}_*(\rmd X, \rmd Y)\,\|X-Y\|^2_{\textsf{Ave}} \right)^{1/2}\\
&=
\inf\limits_{\Pi_* \in \mathcal{C}(\mu,\nu)}\left(\int_{L^2_{\textsf{Ave}} \times L^2_{\textsf{Ave}}} {\Pi}_*(\rmd X, \rmd Y)\,\frac{1}{T}\int_{0}^T \rmd s\, |X(s)-Y(s)|^2 \right)^{1/2}.
\end{split}
\end{equation}
We point out that~\eqref{def:wa} allows us to
connect the paths $X$ and $Y$ on the level of realizations, and 
in fact, this is a natural metric between realizations, see~\cite{BionNadal2019}.
For shorthand, we write 
$\mathcal{W}_2(X,Y)$ in place of $\mathcal{W}_2(\mu, \nu)$. 
The definition given 
in~\eqref{def:wa} defines a metric
that metrizes the weak topology on $\mathcal{P}_T$, see for 
instance~\cite{BionNadal2019}.
For  more details on definitions, properties and notions related to couplings and Wasserstein metrics, we refer to the monographs~\cite{Ambrosio2021},~\cite{Figalli2023},~\cite{Panaretos2020} 
and~\cite{Villani09}.

We start with the following observation. 
For a given observable $F$ 
satisfying~\eqref{eq:wlip} and any coupling $\Gamma\in \mathcal{C}(\mu,\nu)$, the Cauchy--Schwarz inequality and the Minkowski inequality yield for any  $X,Y\in L^2_{\textsf{Ave}}$ the following estimate
\begin{equation}\label{eq:ine1}
\begin{split}
\quad|\mathbb{E}_{\Gamma}[F(X)]-\mathbb{E}_{\Gamma}[F(Y)]|& \leq
\mathbb{E}_{\Gamma}[|F(X)-F(Y)|]\\
&\leq \frac{\kappa}{2}
\mathbb{E}_{\Gamma}[(W(X)+W(Y))\cdot \|X-Y\|_{\textsf{Ave}}]\\
&\leq \frac{\kappa}{2}
\left(\mathbb{E}_{\Gamma}[(W(X)+W(Y))^2]\right)^{1/2}\cdot
\left(\mathbb{E}_{\Gamma}[\|X-Y\|^2_{\textsf{Ave}}]\right)^{1/2}\\
&\leq \frac{\kappa}{2}
((\mathbb{E}_{\Gamma}[(W(X))^2])^{1/2}+(\mathbb{E}_{\Gamma}[(W(Y))^2])^{1/2})\cdot
\left(\mathbb{E}_{\Gamma}[\|X-Y\|^2_{\textsf{Ave}}]\right)^{1/2}, 
\end{split}
\end{equation}
where $\mathbb{E}_{\Gamma}$ denotes the expectation with respect to the probability measure (coupling) $\Gamma$.
Observe that $F(X)$, $F(Y)$, $|W(X)|^2$ and $|W(Y)|^2$ only depend on marginal distributions, hence 
by the definition of coupling, for any coupling $\Gamma\in \mathcal{C}(\mu,\nu)$ we have 
\begin{equation}\label{eq:coudef}
\begin{split}
&\mathbb{E}_{\mu}[F(X)]=\mathbb{E}_{\Gamma}[F(X)],\quad
\mathbb{E}_{\mu}[(W(X))^2]=\mathbb{E}_{\Gamma}[(W(X))^2],\\
&\mathbb{E}_{\nu}[F(Y)]=\mathbb{E}_{\Gamma}[F(Y)],\quad\,\,
\mathbb{E}_{\nu}[(W(Y))^2]=\mathbb{E}_{\Gamma}[(W(Y))^2],
\end{split}
\end{equation}
where $\mathbb{E}_{\mu}$  and $\mathbb{E}_{\nu}$ denote the expectation with respect to the probability measure $\mu$ and $\nu$, respectively.
Then~\eqref{eq:ine1} with the help 
of~\eqref{eq:coudef} yields
\begin{equation}
\begin{split}
|\mathbb{E}_{\mu}[F(X)]-\mathbb{E}_{\nu}[F(Y)]|&\leq \frac{\kappa}{2}
((\mathbb{E}_\mu[(W(X))^2])^{1/2}+(\mathbb{E}_\nu[(W(Y))^2])^{1/2})\cdot\left(\mathbb{E}_{\Gamma}[\|X-Y\|^2_{\textsf{Ave}}]\right)^{1/2}.
\end{split}
\end{equation}
Since $\Gamma\in \mathcal{C}(\mu,\nu)$ is chosen arbitrary,
recalling the definition of $\mathcal{W}_2(X,Y)$ given in~\eqref{def:wa} with its shorthand notation and 
optimizing over all couplings $\Gamma\in \mathcal{C}(\mu,\nu)$ both sides of the inequality~\eqref{eq:ine1}, we obtain the following useful estimate that we state as a proposition.
\begin{proposition}[Weighted Lipschitz observables]\label{lem:ine}
\hfill

\noindent
Let $F$ be an $(\kappa,W)$-Lusin--Lipschitz observable that satisfies~\eqref{eq:wlip}. Then for any $X,Y\in L^2_{\textsf{Ave}}$ it follows that
\begin{equation}\label{eq:estimatecrucial}
\begin{split}
|\mathbb{E}_{\mu}[F(X)]-\mathbb{E}_{\nu}[F(Y)]|\leq \frac{\kappa}{2}
((\mathbb{E}_{\mu}[(W(X))^2])^{1/2}+(\mathbb{E}_{\nu}[(W(Y))^2])^{1/2})\cdot
\mathcal{W}_2(X,Y).
\end{split}
\end{equation}
\end{proposition}

\begin{remark}[Crucial and useful estimate]\hfill

\noindent
The inequality~\eqref{eq:estimatecrucial} allows us to estimate the difference of the moments for Lusin--Lipschitz observables at paths $X$ and $Y$ by the Wasserstein distance between $X$ and $Y$, and the average of weight functions of the individual paths $X$ and $Y$. 
Recently, using probabilistic methods such as couplings (coin-flip coupling, synchronous coupling and Koml\'os--Major--Tusn\'ady coupling), a non-asymptotic estimation of $\mathcal{W}_2(X,Y)$ is provided in~\cite{BarreraLukkarinen2023}. 

Let $\Gamma\in \mathcal{C}(\mu,\nu)$ be fixed and $F$ be an $(\kappa,W)$-Lusin--Lipschitz observable that satisfies~\eqref{eq:wlip}.
Inspecting the inequality~\eqref{eq:ine1} 
with the help of the H\"older inequality and the Minkowski inequality
one can deduce that
\begin{equation}\label{eq:inenueva}
\begin{split}
\quad|\mathbb{E}_{\mu}[F(X)]-\mathbb{E}_{\nu}[F(Y)]|& \leq \frac{\kappa}{2}
((\mathbb{E}_{\mu}[(W(X))^q])^{1/q}+(\mathbb{E}_{\nu}[(W(Y))^q])^{1/q})\cdot
\left(\mathbb{E}_{\Gamma}[\|X-Y\|^p_{\textsf{Ave}}]\right)^{1/p}, 
\end{split}
\end{equation}
for $p>1$ and $q>1$ such that $1/p+1/q=1$, whenever the 
right-hand side of~\eqref{eq:inenueva} exists.
In particular, for $p=q=2$ we obtain Proposition~\ref{lem:ine}.

\end{remark}

\subsection{\textbf{The main results and its consequences}}\label{sub:results}
\hfill

\noindent
In this subsection, we state the main results
of this manuscript and its consequences.
Recall that  $X^{\textsf{sym},v}:=(X^{\textsf{sym},v}(t))_{t\geq 0}$ denotes the symmetric telegraph process with velocity parameter $v\neq 0$ and rate $\lambda>0$.

The first main result of this manuscript is the following.
\begin{theorem}[Symmetric Goldstein--Kac telegraph process]\label{thA}
\hfill

\noindent
Assume that $\lambda=\lambda^*>0$ and $v_0=-v^*_0$ with $v_0\neq 0$.
Let $F$ be a $(\kappa,W)$-Lipschitz observable that satisfies~\eqref{eq:wlip}.  
Then there exists a positive (absolute) constant $C$ such that for any $L>0$, $T>0$, $\lambda>0$, $v_0\neq 0$ it follows that
\begin{equation}\label{eq:inec1}
\begin{split}
|\mathbb{E}_{\mu}[F(L^{-1}X^{\textsf{sym},v_0})]-\mathbb{E}_{\nu}[F(B)]|&\leq \frac{\kappa}{2}
\mathcal{W}_2(L^{-1}X^{\textsf{sym},v_0},B)\\
&\qquad \times
((\mathbb{E}_{\mu}[(W(L^{-1}X^{\textsf{sym},v_0}))^2])^{1/2}+(\mathbb{E}_{\nu}[(W(B))^2])^{1/2})
\end{split}
\end{equation}
with
\begin{equation}\label{eq:mainresult}
\mathcal{W}_2(L^{-1}X^{\textsf{sym},v_0},B)\leq C\sqrt{{T_*}{L^{-2}_*}}T^{-1/4}_*\big(\sqrt{\ln(T_*+3)}+T^{-3/4}_*\big)+CL^{-1}_*,
\end{equation}
where the constants 
$T_*$ and $L_*$ are given by
\begin{equation}\label{e:scalings2020}
T_*:=\lambda T,\qquad {L}_*:=|v_0|^{-1}\lambda L,
\end{equation}
and the diffusivity constant of the Brownian motion (in a suitable probability space) $B:=(B(t))_{t\geq 0}$ is defined by 
\begin{equation}\label{eq:diffusivity}
\sigma^2:=L^{-2}\frac{v_0^2}{\lambda}=L^{-2}_*\lambda\,.
\end{equation}
\end{theorem}
The proof is presented in Subsection~\ref{sec:pthA}.

In the sequel, we study the analogous of Theorem~\ref{thA} for asymmetric telegraph processes and state it as a theorem.

\begin{theorem}[Asymmetric Goldstein--Kac telegraph process: same rate]\label{th:asysamerate}
\hfill

\noindent
Assume that $\lambda=\lambda^*>0$ and $v_0+v^*_0\neq 0$, and set $v:=\frac{v_0+v^*_0}{2}$.
Let $F$ be a $(\kappa,W)$-Lipschitz observable that satisfies~\eqref{eq:wlip}.
Then there exists
a positive absolute constant $C$ such that for any $L>0$, $T>0$, $\lambda>0$, $v\neq 0$ it holds
\begin{equation}
\begin{split}
|\mathbb{E}_{\mu}[F(L^{-1}X)]-\mathbb{E}_{\nu}[F(\widetilde{B})]|&\leq \frac{\kappa}{2}
\mathcal{W}_2(L^{-1}X,\widetilde{B})\\
&\qquad \times
((\mathbb{E}_{\mu}[(W(L^{-1}X)^2])^{1/2}+(\mathbb{E}_{\nu}[(W(\widetilde{B}))^2])^{1/2}),
\end{split}
\end{equation}
where  $\widetilde{B}$ is a Brownian motion with drift $\textsf{d}:=\frac{v_0-v^*_0}{2L}$  and diffusivity constant $\sigma^2$ defined in~\eqref{eq:diffusivity}. Moreover, 
\begin{equation}
\mathcal{W}_2(L^{-1}X,\widetilde{B})
\leq C\sqrt{{\widetilde{T}_*}{\widetilde{L}^{-2}_*}}\widetilde{T}^{-1/4}_*\Big(\sqrt{\ln(\widetilde{T}_*+3)}+\widetilde{T}^{-3/4}_*\Big)+C\widetilde{L}^{-1}_*,
\end{equation}
where the constants 
$\widetilde{T}_*$ and $\widetilde{L}_*$ are given by
\begin{equation}
\widetilde{T}_*:=\lambda T,\qquad \widetilde{L}_*:=|v|^{-1}\lambda L.
\end{equation}
\end{theorem}

The proof is given in Subsection~\ref{sec:pasy}.

Now, we  provide a class of exponential functionals that satisfy~\eqref{eq:wlip}. 
We stress that the parameter $\kappa$
and the weight function $W$ are given explicitly.
\begin{lemma}[Exponential weighted Lipschitz observables]\label{lemF}
\hfill

\noindent
Let $f:\mathbb{R}\to \mathbb{R}$  and $g:\mathbb{R}\to \mathbb{R}$  be Lipschitz functions with some Lipschitz constants $\kappa_f$ and $\kappa_g$, respectively.
For a given constant $a\in \mathbb{R}$ consider the set of functions
\begin{equation}\label{eq:2aX}
U_a:=\{X\in L^2_{\textsf{Ave},\mathbb{R}}: e^{aX}\in L^2_{\textsf{Ave},\mathbb{R}}\}
\end{equation}
and then for $b\in \mathbb{R}$ define
 $F_{a,b}:U_a\to \mathbb{R}$ by 
\begin{equation}\label{defFab}
F_{a,b}(X):=g\left(\frac{1}{T}\int_{0}^T \rmd s \,f\big(e^{aX(s)+bs}\big)\right).
\end{equation}
Then $F_{a,b}$ satisfies~\eqref{eq:wlip} with 
\begin{equation}
\kappa:=2|a|\kappa_f\kappa_g\quad \textrm{ and } \quad W(X)=W_{a,b}(X):=\sqrt{\frac{1}{T}\int_{0}^T \rmd s \,e^{2aX(s)+2bs}}.
\end{equation}
In particular, for a given constants $\ell,c\in \mathbb{R}$ the functions $f(z):=\ell z$, $z\in \mathbb{R}$ and  $f(z):=\max\{\ell z,c\}$, $z\in \mathbb{R}$ may have Lipschitz constants  $\kappa_f=\ell$.
\end{lemma}
The proof is provided in Subsection~\ref{sec:plemF}.

\begin{remark}[More general observables] 
\hfill

\noindent
Following step by step the proof of Lemma~\ref{lemF}, one can see that the observable
\begin{equation}\label{defFab1}
F(X):=g\left(\frac{1}{T}\int_{0}^T \rmd s \,f\big(e^{aX(s)}h(s)\big)\right)
\end{equation}
with $h:\mathbb{R}\to \mathbb{R}$ being a continuous function
satisfying~\eqref{eq:wlip} with 
\begin{equation}\label{eq:kappaWfor}
\kappa:=2|a|\kappa_f\kappa_g\quad \textrm{ and } \quad W(X):=\sqrt{\frac{1}{T}\int_{0}^T \rmd s \,e^{2aX(s)}(h(s))^2}.
\end{equation}
\end{remark}

Generically, the precise formula for the distribution of $(W(X))^2$ given in~\eqref{eq:kappaWfor} is  hard to obtain,
see~\cite{Kac1949}. Nevertheless, we only require its expectation, which is an easier task after applying Fubini's Theorem.
As a consequence of Theorem~\ref{thA} and Lemma~\ref{lemF} we have the following estimates for path-dependent exponential weighted Lipschitz observables.
\begin{corollary}[Error estimates: symmetric motion]\label{corA}
\hfill

\noindent
Suppose that all assumptions and notation made in Lemma~\ref{lemF} hold. 
Assume that $\lambda=\lambda^*>0$ and $v_0=-v^*_0$ with $v_0\neq 0$.
Then there exists a positive absolute constant $C$ such that for any $L>0$, $T>0$, $\lambda>0$, $v_0\neq 0$ it follows that
\begin{equation}
\begin{split}
\left|\mathbb{E}_{\mu}[F_{a,b}(L^{-1}X^{\textsf{sym},v_0})]-\mathbb{E}_{\nu}[F_{a,b}(B)]\right|&\leq |a|\kappa_f\kappa_g
\mathcal{W}_2(L^{-1}X^{\textsf{sym},v_0},B)\\
&\hspace{-0.5cm}\times
((\mathbb{E}_{\mu}[(W_{a,b}(L^{-1}X^{\textsf{sym},v_0}))^2])^{1/2}+(\mathbb{E}_{\nu}[(W_{a,b}(B))^2])^{1/2})
\end{split}
\end{equation}
where
\begin{equation}
\mathcal{W}_2(L^{-1}X^{\textsf{sym},v_0},B)\leq C\sqrt{{T_*}{L^{-2}_*}}T^{-1/4}_*\big(\sqrt{\ln(T_*+3)}+T^{-3/4}_*\big)+CL^{-1}_*,
\end{equation}
\begin{equation}\label{eq:WabB}
\mathbb{E}_{\nu}\left[(W_{a,b}(B))^2\right]=\frac{1}{2a^2 T_*L^{-2}_*+2bT_*/\lambda}\left(e^{2a^2 T_*L^{-2}_*+2bT_*/\lambda}-1\right),
\end{equation}
\begin{equation}\label{eq:WabX}
\begin{split}
\mathbb{E}_{\mu}[(W_{a,b}(L^{-1}X^{\textsf{sym},v_0}))^2]&\leq \left(1+\frac{1}{\sqrt{1+4a^2L^{-2}_*}}\right)\\
&\quad\times\max\left\{1,\exp\left({2T_*b/\lambda+\frac{4a^2T_*L^{-2}_*}{1+\sqrt{1+4a^2L^{-2}_*}}}\right)\right\}.
\end{split}
\end{equation}
Here, the constants 
$T_*$ and $L_*$ are given by
$T_*:=\lambda T,$ ${L}_*=|v_0|^{-1}\lambda L$,
and the diffusivity constant of the Brownian motion $B:=(B(t))_{t\geq 0}$ (in a suitable probability space) is defined by 
$\sigma^2:=L^{-2}\frac{v_0^2}{\lambda}$.

The natural choice  $b=-\frac{1}{2}a^2\sigma^2$ as in the risk-neutral measure in the Black--Scholes--Merton model yields 
\begin{equation}
\begin{split}
\mathbb{E}_{\nu}\left[(W_{a,b}(B))^2\right]&=\frac{1}{a^2 T_*L^{-2}_*}\left(e^{a^2 T_*L^{-2}_*}-1\right),\\
\mathbb{E}_{\mu}[(W_{a,b}(L^{-1}X^{\textsf{sym},v_0}))^2]
&\leq \left(1+\frac{1}{\sqrt{1+4a^2L^{-2}_*}}\right)\\
&\quad \times \max\left\{1,\exp\left(a^2T_*L^{-2}_*\left(\frac{4-\sqrt{1+4a^2L^{-2}_*}}{\sqrt{1+4a^2L^{-2}_*}}\right)\right)\right\}.
\end{split}
\end{equation}
In particular,
\begin{equation}
\mathbb{E}_{\mu}[(W_{a,b}(L^{-1}X^{\textsf{sym},v_0}))^2]
\leq 2\max\left\{1,\exp\left(4a^2T_*L^{-2}_*\right)\right\}.
\end{equation}
\end{corollary}

The proof is presented in Subsection~\ref{sec:pcorA}. We point out that the computation of the left-hand side of~\eqref{eq:WabX} requires the shape of the spatial Laplace transform (a.k.a. Moment generating function) for the symmetric Goldstein--Kac telegraph process, which is given in terms of 
the hyperbolic cosine and hyperbolic sine functions,
see Theorem~4 in~\cite{Kolesnik2012}.

Now, we state the analogous of Corollary~\ref{corA} when $\lambda=\lambda^*$ and $v_0+v^*_0\neq 0$.

\begin{corollary}[Error estimates: asymmetric motion]\label{corB}
\hfill

\noindent
Suppose that all assumptions and notation made in Lemma~\ref{lemF} hold. 
Assume that $\lambda=\lambda^*>0$ and $v_0+v^*_0\neq 0$, and set $v:=\frac{v_0+v^*_0}{2}$.
Then there exists
a positive absolute constant $C$ such that for any $L>0$, $T>0$, $\lambda>0$, $v\neq 0$ it holds
\begin{equation}
\begin{split}
|\mathbb{E}_{\mu}[F_{a,b}(L^{-1}X)]-\mathbb{E}_{\nu}[F_{a,b}(\widetilde{B})]|&\leq \sqrt{2}|a|\kappa_f \kappa_g
\mathcal{W}_2(L^{-1}X,\widetilde{B})\\
&\qquad \times
((\mathbb{E}_{\mu}[(W_{a,b}(L^{-1}X)^2])^{1/2}+(\mathbb{E}_{\nu}[(W_{a,b}(\widetilde{B}))^2])^{1/2}),
\end{split}
\end{equation}
where $\widetilde{B}$ is a Brownian motion with drift $\textsf{d}:=\frac{v_0-v^*_0}{2L}$  and the diffusivity constant $\sigma^2$ is defined 
in~\eqref{eq:diffusivity}.
Moreover,
\begin{equation}
\mathcal{W}_2(L^{-1}X,\widetilde{B})
\leq C\sqrt{{\widetilde{T}_*}{\widetilde{L}^{-2}_*}}\widetilde{T}^{-1/4}_*\Big(\sqrt{\ln(\widetilde{T}_*+3)}+\widetilde{T}^{-3/4}_*\Big)+C\widetilde{L}^{-1}_*,
\end{equation}
where $\widetilde{T}_*:=\lambda T$, $\widetilde{L}_*:=|v|^{-1}\lambda L$.
In addition,
\begin{equation}\label{eq:cambio}
\begin{split}
\mathbb{E}_{\nu}\left[(W_{a,b}(\widetilde{B}))^2\right]&=\mathbb{E}\left[(W_{a,b+\textsf{d}}(B))^2\right],\\
\mathbb{E}_{\mu}[(W_{a,b}(L^{-1}X))^2]
&=\mathbb{E}[(W_{a,b+\textsf{d}}(L^{-1}X^{\textsf{sym},v}))^2],
\end{split}
\end{equation}
where $B:=(B(t))_{t\geq 0}$ is a Brownian motion with diffusivity constant  
$\sigma^2$, and the right-hand side 
of~\eqref{eq:cambio} can be estimated 
by~\eqref{eq:WabB} and~\eqref{eq:WabX}.
\end{corollary}

The proof is given in Subsection~\ref{sub:pcorB}.

\begin{remark}[Weighted Lipschitz observables that are not integrable]
\hfill

\noindent
There are weighted Lipschitz observables $F$ which are not integrable and hence Corollary~\ref{corA} is not meaningful.
For instance,
for any $a>0$ , it is not hard to see that $F:L^2_{\textsf{Ave},\mathbb{R}}\to [0,\infty)$ given by 
\[
F(X):=e^{a T^{-1}\int_{0}^T \rmd s\, (X(s))^2}=e^{a\|X\|^2_{\textsf{Ave}}}
\]
satisfies~\eqref{eq:wlip} with
\[
\kappa:=4a \quad \textrm{ and }\quad W(X):=e^{a\|X\|^2_{\textsf{Ave}}}\|X\|_{\textsf{Ave}}.
\]
We note that the crude estimate $|X^{\textsf{sym}}(s)|\leq |v_0|s$ for any $t\geq 0$ implies that the expectation $\mathbb{E}_{\mu}[(W(L^{-1}X^{\textsf{sym}}))^2]$ always exists.
However, 
\begin{equation}\label{eq:condiciones}
\begin{split}
&\textrm{ for all }\quad a>\frac{3}{2}\frac{1}{L^2_*T_*}\quad \textrm{ it follows that }\quad\mathbb{E}_{\nu}[W^2(B)]=\infty,\\
&
\textrm{ whereas }\quad \mathbb{E}_{\nu}[W^2(B)]<\infty\quad \textrm{ for }\quad 0<a<\frac{1}{2}\frac{1}{L^2_*T_*}.
\end{split}
\end{equation} 
In the sequel, we show that  $\mathbb{E}_{\nu}[W^2(B)]<\infty$ for $0<a<\frac{1}{2}\frac{1}{L^2_*T_*}$. Indeed, the Monotone Convergence Theorem yields
\begin{equation}
\begin{split}
\mathbb{E}_{\nu}[F(B)]&=\mathbb{E}_{\nu}\left[e^{a T^{-1}\int_{0}^T \rmd s\, (B(s))^2}\right]=\sum_{n=0}^{\infty}\frac{a^n}{n!}\mathbb{E}_{\nu}\left[\left(T^{-1}\int_{0}^T \rmd s\, (B(s))^2\right)^n\right].
\end{split}
\end{equation}
Now, we apply the Jensen inequality for each $n\in \mathbb{N}_0$ and obtain
\begin{equation}
\begin{split}
\left(T^{-1}\int_{0}^T \rmd s\, (B(s))^2\right)^n\leq 
\left(T^{-1}\int_{0}^T \rmd s\, (B(s))^{2n}\right),
\end{split}
\end{equation}
which with the help of the Fubini Theorem implies
\begin{equation}
\begin{split}
\sum_{n=0}^{\infty}\frac{a^n}{n!}\mathbb{E}_{\nu}\left[\left(T^{-1}\int_{0}^T \rmd s\, (B(s))^2\right)^n\right]&\leq  \sum_{n=0}^{\infty}\frac{a^n}{n!}\mathbb{E}_{\nu}\left[T^{-1}\int_{0}^T \rmd s\, (B(s))^{2n}\right]\\
&=\sum_{n=0}^{\infty}\frac{a^n}{n!}T^{-1}\int_{0}^T \rmd s\, \mathbb{E}_{\nu}\left[(B(s))^{2n}\right]\\
&=\sum_{n=0}^{\infty}\frac{a^n}{n!}\sigma^{2n} \mathbb{E}[\mathsf{N}^{2n}]T^{-1}\int_{0}^T \rmd s\,  s^n \\
&=\sum_{n=0}^{\infty}\frac{a^n}{(n+1)!}\sigma^{2n} \mathbb{E}[\mathsf{N}^{2n}]T^n,
\end{split}
\end{equation}
where $\mathsf{N}$ denotes the standard Gaussian distribution.
Since 
\[
\mathbb{E}[\mathsf{N}^{2n}]=\prod_{j=1}^{n}(2j-1)=:(2n-1)!!\quad \textrm{ and }\quad
\sigma^2T=L^{-2}\frac{v_0^2}{\lambda^2}\lambda T=L^{-2}_*T_*,
\]
we have
\begin{equation}\label{eq:series}
\begin{split}
\sum_{n=0}^{\infty}\frac{a^n}{n!}\mathbb{E}_{\nu}\left[\left(T^{-1}\int_{0}^T \rmd s\, (B(s))^2\right)^n\right]&\leq 1+ \sum_{n=1}^{\infty}\frac{(aL^{-2}_*T_*)^n}{(n+1)!} (2n-1)!!,
\end{split}
\end{equation}
where $k!!$ denotes the double factorial of a given natural number $k$.
Using the celebrated ratio test for series we deduce that the series in the right-hand side of~\eqref{eq:series} converges absolutely when $0<a<\frac{1}{2}\frac{1}{L^2_*T_*}$.

We continue with the proof of 
$\mathbb{E}_{\nu}[W^2(B)]=\infty$ when $a>\frac{3}{2}\frac{1}{L^2_*T_*}$.
We recall that
\begin{equation}\label{eq:BB}
W_T(s):=B(s)-\frac{s}{T}B(T), \quad s\in [0,T]
\end{equation}
defines a Brownian bridge on $[0,T]$, that is, $W_T(0)=W_T(T)=0$.
It is known that $B(T)$ is independent of $(W_T(s))_{s\in [0,T]}$, see for instance Subsection~1.1 in~\cite{Mansuy2008}.
We note that 
\begin{equation}
\begin{split}
e^{a T^{-1}\int_{0}^T \rmd s\, (B(s))^2}&=e^{a T^{-1}\int_{0}^T \rmd s\, (W_T(s)+B(T)T^{-1}s)^2}\\
&=
e^{a T^{-1}\int_{0}^T \rmd s\, (W_T(s))^2 }
e^{a T^{-1}\int_{0}^T \rmd s\, (B(T)T^{-1}s)^2}
e^{a T^{-1}\int_{0}^T \rmd s\, 2W_T(s)B(T)T^{-1}s}\\
&=
e^{a T^{-1}\int_{0}^T \rmd s\, (W_T(s))^2 }
e^{a (B(T))^2 /3}
e^{a T^{-2}B(T)\int_{0}^T \rmd s\, 2W_T(s)s}.
\end{split}
\end{equation}
Then taking expectation in both sides and using the 
law of total expectation
(a.k.a. tower rule, see Theorem~8.14 of~\cite{Klenke2020}) of conditional expectation, the Jensen inequality for conditional expectation (see~Theorem~8.20 in~\cite{Klenke2020}), the fact that  $B(T)$ is independent of $(W_T(s))_{s\in [0,T]}$ and the Fubini Theorem we obtain
\begin{equation}\label{eq:inec}
\begin{split}
\mathbb{E}\left[e^{a T^{-1}\int_{0}^T \rmd s\, (B(s))^2}\right]
&=\mathbb{E}\left[e^{a T^{-1}\int_{0}^T \rmd s\, (W_T(s))^2 }
e^{a (B(T))^2 /3}
e^{a T^{-2}B(T)\int_{0}^T \rmd s\, 2W_T(s)s}\right]\\
&=\mathbb{E}\left[\mathbb{E}\left[e^{a T^{-1}\int_{0}^T \rmd s\, (W_T(s))^2 }
e^{a (B(T))^2 /3}
e^{a T^{-2}B(T)\int_{0}^T \rmd s\, 2W_T(s)s}\Big|B(T)\right]\right]\\
&\geq 
\mathbb{E}\left[e^{a (B(T))^2 /3}
\mathbb{E}\left[
e^{a T^{-2}B(T)\int_{0}^T \rmd s\, 2W_T(s)s}\Big|B(T)\right]\right]\\
&\geq 
\mathbb{E}\left[e^{a (B(T))^2 /3}
e^{\mathbb{E}\left[a T^{-2}B(T)\int_{0}^T \rmd s\, 2W_T(s)s\Big|B(T)\right]}\right]\\
&=\mathbb{E}\left[e^{a (B(T))^2 /3}
e^{a T^{-2}B(T)\mathbb{E}\left[\int_{0}^T \rmd s\, 2W_T(s)s\Big|B(T)\right]}\right]\\
&=\mathbb{E}\left[e^{a (B(T))^2 /3}
e^{a T^{-2}B(T)\mathbb{E}\left[\int_{0}^T \rmd s\, 2W_T(s)s\right]}\right]\\
&=\mathbb{E}\left[e^{a (B(T))^2 /3}\right].
\end{split}
\end{equation}
Now, we compute the right-hand side of 
the~\eqref{eq:inec}, that is,
\begin{equation}
\mathbb{E}\left[e^{a (B(T))^2 /3}\right]=
\mathbb{E}\left[e^{a T\sigma^2 \mathsf{N}^2/3}\right],
\end{equation}
where $\mathsf{N}$ denotes the standard Gaussian distribution. Then we have
\begin{equation}
\begin{split}
\mathbb{E}\left[e^{a T\sigma^2 \mathsf{N}^2/3}\right]=
\frac{1}{\sqrt{2\pi}}\int_{\mathbb{R}}
e^{a T\sigma^2 x^2/3}e^{-x^2/2}\rmd x=
\frac{1}{\sqrt{2\pi}}\int_{\mathbb{R}}
e^{x^2(a T\sigma^2/3-1/2)}\rmd x=\infty
\end{split}
\end{equation}
when $a T\sigma^2/3-1/2>0$. In other words, when $aL^{-2}_*T_*>3/2$.

For any non-negative   $a$, the Monotone Convergence Theorem implies 
\[
\mathbb{E}_{\nu}\left[e^{a T^{-1}\int_{0}^T \rmd s\, (B(s))^2}\right]=\lim\limits_{N\to \infty}\sum_{n=0}^{N} \frac{a^n}{n!}c_n,
\]
where 
\[
c_n:=\mathbb{E}_{\nu}\left[\left(T^{-1}\int_{0}^T \rmd s\, (B(s))^2\right)^n\right], \quad n\in \mathbb{N}\cup\{0\}.
\]
By~\eqref{eq:series} we have that the function
\[
a\mapsto \psi(a):=\mathbb{E}_{\nu}\left[e^{a T^{-1}\int_{0}^T \rmd s\, (B(s))^2}\right]
\]
is analytic for $a\in (-\infty,\frac{1}{2}\frac{1}{L^2_*T_*})$. 
Therefore, there exists $r>0$ such that 
\[
\sum_{n=0}^{\infty} \frac{a^n}{n!}c_n <\infty\quad \textrm{ whenever }\quad |a|<r.
\]
By~\eqref{eq:condiciones} we have that $\frac{1}{2}\frac{1}{L^2_*T_*}\leq r<\frac{3}{2}\frac{1}{L^2_*T_*}$.
\end{remark}

The rest of the manuscript is organized as follows. 
In Section~\ref{sec:proofs}
we present the proofs of the results given in Section~\ref{sec:intr}. More precisely, in Subsection~\ref{sec:pthA} we show Theorem~\ref{thA}, in Subsection~
\ref{sec:pasy} we prove
Theorem~\ref{th:asysamerate}, in Subsection~\ref{sec:plemF} we show
Lemma~\ref{lemF}, in Subsection~\ref{sec:pcorA}
we prove Corollary~\ref{corA}, and finally in Subsection~\ref{sub:pcorB} we give the steps for the proof of Corollary~\ref{corB}.
In Section~\ref{sec:simulations} we present results of simulations, and assess the sharpness of the error bounds.

\section{\textbf{Proofs of the results}}\label{sec:proofs}
In this section, we present the proofs of the results stated in Subsection~\ref{sub:results}.

\subsection{\textbf{Proof of Theorem~\ref{thA}}}\label{sec:pthA}\hfill

\noindent
Since $F$ is a $(\kappa,W)$-Lipschitz observable that satisfies~\eqref{eq:wlip}, Proposition~\ref{lem:ine} implies~\eqref{eq:inec1}. While~
\eqref{eq:mainresult} is a consequence of Corollary~2.7 
in~\cite{BarreraLukkarinen2023}.
\qedw

\subsection{\textbf{Proof of Theorem~\ref{th:asysamerate}}}\label{sec:pasy}\hfill

\noindent
Using the so-called Galilean transformation it is possible  to symmetrize the motion. To be more precise,
by Formula~(1.1) and~(1.2) 
in~\cite{Cinque2021} for the asymmetric telegraph process $(X(t))_{t\geq 0}$ can be written as 
\begin{equation}
\begin{split}
X(t)&=\left(\frac{v_0-v^*_0}{2}\right)t+\left(V(0)-\frac{v_0-v^*_0}{2}\right)\int_0^t \rmd s\, (-1)^{N(s)}\\
&=\left(\frac{v_0-v^*_0}{2}\right)t\\
&\qquad+\left(v_0\,\mathsf{1}(V(0)=v_0)
+(-v^*_0)\,\mathsf{1}(V(0)=-v^*_0)-\frac{v_0-v^*_0}{2}\right)\int_0^t \rmd s\, (-1)^{N(s)}\\
&=\left(\frac{v_0-v^*_0}{2}\right)t+\left(\frac{v_0+v^*_0}{2}\right)\left(\mathsf{1}(V(0)=v_0)-\mathsf{1}(V(0)=-v^*_0)\right)\int_0^t \rmd s\, (-1)^{N(s)}\\
&=\left(\frac{v_0-v^*_0}{2}\right)t+\left(\mathsf{1}(V(0)=v_0)-\mathsf{1}(V(0)=-v^*_0)\right)X^{\textsf{sym},v}(t)\quad \textrm{ with }\quad v:=\frac{v_0+v^*_0}{2},
\end{split}
\end{equation}
where $N:=(N(t))_{t\geq 0}$ is Poisson process with positive intensity $\lambda$, $V(0)$ is a random variable taking values in the state space $\{v_0,-v^*_0\}$ with uniform probability $\mu$ (i.e., $\mu(\{v_0\})=\mu(\{-v^*_0\})=1/2$) and independent of the process $N$.

We note that $|\mathsf{1}(V(0)=v_0)-\mathsf{1}(V(0)=-v^*_0)|=1$ and recall $B$ and $-B$ have the same distribution, where $B$ is given in Theorem~\ref{thA}.
Let $\textsf{d}=\frac{1}{L}\frac{v_0-v^*_0}{2}$ and define $\widetilde{B}$ being the Brownian motion with drift $a$ and diffusivity constant $\sigma^2$ given 
in~\eqref{eq:diffusivity}. In other words,
\[
\widetilde{B}(t)\stackrel{\mathrm{Law}}{=}\textsf{d}t+B(t),\quad t\geq 0,
\] 
where 
$\stackrel{\mathrm{Law}}{=}$ denotes equality in distribution and
$B$ is the Brownian motion defined in Theorem~\ref{thA}.
Proposition~\ref{lem:ine} 
yields  
the existence of a positive (absolute) constant $C$ such that for any $L>0$, $T>0$, $\lambda>0$, $|v|>0$ we have
\begin{equation}
\begin{split}
|\mathbb{E}_{\mu}[F(L^{-1}X)]-\mathbb{E}_{\nu}[F(\widetilde{B})]|&\leq \frac{\kappa}{\sqrt{2}}
\mathcal{W}_2(L^{-1}X,\widetilde{B})\cdot
((\mathbb{E}_{\mu}[(W(L^{-1}X)^2])^{1/2}+(\mathbb{E}_{\nu}[(W(\widetilde{B}))^2])^{1/2}).
\end{split}
\end{equation}
In addition, 
Corollary~2.7 
in~\cite{BarreraLukkarinen2023} implies
\begin{equation}
\mathcal{W}_2(L^{-1}X,\widetilde{B})=\mathcal{W}_2(L^{-1}X^{\textsf{sym},v},B)
\leq C\sqrt{{\widetilde{T}_*}{\widetilde{L}^{-2}_*}}\widetilde{T}^{-1/4}_*\Big(\sqrt{\ln(\widetilde{T}_*+3)}+\widetilde{T}^{-3/4}_*\Big)+C\widetilde{L}^{-1}_*,
\end{equation}
where the constants 
$\widetilde{T}_*$ and $\widetilde{L}_*$ are given by
\begin{equation}
\widetilde{T}_*:=\lambda T,\qquad \widetilde{L}_*=|v|^{-1}\lambda L.
\end{equation}

\subsection{\textbf{Proof of Lemma~\ref{lemF}}}
\label{sec:plemF}
\hfill

\noindent
Note that~\eqref{eq:2aX} implies that $e^{aX}\in L^1_{\textsf{Ave},\mathbb{R}}$.
Since $f$ and $g$ are Lipschitz functions, we have that $F$ is well-defined.
For any $X,Y\in U_a$, the Cauchy--Schwarz  inequality yields
\begin{equation}\label{eq:u1}
\begin{split}
|F(X)-F(Y)|& =\left|g\left(\frac{1}{T}\int_{0}^T \rmd s \,f\big(e^{aX(s)+bs}\big)\right)-g\left(\frac{1}{T}\int_{0}^T \rmd s \,f\big(e^{aY(s)+bs}\big)\right)\right|\\
& \leq \kappa_g
\left|\frac{1}{T}\int_{0}^T \rmd s \,f\big(e^{aX(s)+bs}\big)-\frac{1}{T}\int_{0}^T \rmd s \,f\big(e^{aY(s)+bs}\big)\right|\\
& \leq \kappa_g
\frac{1}{T}\int_{0}^T \rmd s \, \left|f\big(e^{aX(s)+bs}\big)-f\big(e^{aY(s)+bs}\big)\right|\\
&\leq \kappa_g
\frac{1}{T}\int_{0}^T \rmd s \,\kappa_{f}\,|e^{aX(s)+bs}-e^{aY(s)+bs}|.
\end{split}
\end{equation}
We claim that
\begin{equation}\label{eq:expine}
|e^{ax}-e^{ay}|\leq \max\{e^{ax},e^{ay}\}|ax-ay|\leq |a|(e^{ax}+e^{ay})|x-y|\quad \textrm{ for any }\quad x,y\in \mathbb{R}.
\end{equation}
Indeed, let $x,y\in \mathbb{R}$ and assume that $ay\leq ax$. Recall that $1-e^{-z}\leq z$ for all $z\in \mathbb{R}$. Then we have 
\begin{equation}
\begin{split}
|e^{ax}-e^{ay}|&=
e^{ax}|1-e^{-(ax-ay)}|
=e^{ax}(1-e^{-(ax-ay)})
\leq e^{ax}(ax-ay)\\
&\leq 
e^{ax}|ax-ay|\leq 
|a|(e^{ax}+e^{ay})|x-y|,
\end{split}
\end{equation}
while the case $ax\leq ay$ follows similarly.

By~\eqref{eq:expine} we have
\begin{equation}\label{eq:u2}
\begin{split}
\frac{1}{T}\int_{0}^T & \rmd s \,|e^{aX(s)+bs}-e^{aY(s)+bs}| =\frac{1}{T}\int_{0}^T \rmd s \,e^{bs}|e^{aX(s)}-e^{aY(s)}|\\
&\leq |a|\frac{1}{T}\int_{0}^T \rmd s \,e^{bs}(e^{aX(s)}+e^{aY(s)})|X(s)-Y(s)|\\
&= |a|\frac{1}{T}\int_{0}^T \rmd s \,e^{bs}e^{aX(s)}|X(s)-Y(s)|
+|a|\frac{1}{T}\int_{0}^T \rmd s \,e^{bs}e^{aY(s)}|X(s)-Y(s)|\\
&\leq |a|\sqrt{\frac{1}{T}\int_{0}^T \rmd s \,e^{2bs}e^{2aX(s)}}\, \|X-Y\|_{\textsf{Ave}}+ |a|\sqrt{\frac{1}{T}\int_{0}^T \rmd s \,e^{2bs}e^{2aY(s)}}\, \|X-Y\|_{\textsf{Ave}}\\
&=|a|(W(X)+W(Y))\|X-Y\|_{\textsf{Ave}},
\end{split}
\end{equation}
where in the last inequality we use the Cauchy--Schwarz  inequality. 
By~\eqref{eq:u1}  and~\eqref{eq:u2} we deduce the statement.
\qedw

\subsection{\textbf{Proof of Corollary~\ref{corA}}}
\label{sec:pcorA}
\hfill

\noindent
For $a=0$, the statement holds true trivially. Hence we always assume that $a\neq 0$.
By Lemma~\ref{lemF} we only need to estimate (an upper bound) of $\mathbb{E}_{\mu}[(W_{a,b}(L^{-1}X^{\textsf{sym},v_0}))^2]$ and $\mathbb{E}_{\mu}[(W_{a,b}(B)^2]$.

We start with the computation of  $\mathbb{E}_{\mu}[(W_{a,b}(B)^2]$.
By the Fubini Theorem we have
\begin{equation}
\begin{split}
\mathbb{E}_{\nu}\left[(W_{a,b}(B))^2\right]&=
\mathbb{E}_{\nu}\left[T^{-1}\int_{0}^T \rmd s \,e^{2aB(s)+2bs}\right] =
T^{-1}\int_{0}^T \rmd s \,\mathbb{E}_{\nu}\left[e^{2aB(s)+2bs}\right]\\
&=T^{-1}\int_{0}^T \rmd s \,e^{2bs}\mathbb{E}_{\nu}\left[e^{2aB(s)}\right]=T^{-1}\int_{0}^T \rmd s \,e^{2bs}\mathbb{E}\left[e^{2a\sigma \sqrt{s}\mathsf{N}}\right]\\
&=
T^{-1}\int_{0}^T \rmd s \,e^{2bs}e^{2a^2\sigma^2 s}=\frac{1}{2T(a^2\sigma^2+b)}\left(e^{2(a^2\sigma^2+b)T}-1\right).
\end{split}
\end{equation}
where $\mathsf{N}$ denotes a standard Gaussian distribution.

We continue with the estimation of $\mathbb{E}_{\mu}[(W_{a,b}(L^{-1}X^{\textsf{sym},v_0}))^2]$. Again, 
by the Fubini Theorem we have
\begin{equation}\label{eq:fb}
\begin{split}
\mathbb{E}_{\mu}[(W_{a,b}(L^{-1}X^{\textsf{sym},v_0}))^2]&=
\mathbb{E}_{\mu}\left[T^{-1}\int_{0}^T \rmd s \,e^{2aL^{-1}X^{\textsf{sym},v_0}(s)+2bs}\right]\\
&=
T^{-1}\int_{0}^T \rmd s \,e^{2bs}\mathbb{E}_{\mu}\left[e^{2aL^{-1}X^{\textsf{sym},v_0}(s)}\right].
\end{split}
\end{equation}
Now, doing the change of variable $z\mapsto z^2$ in 
Theorem~4 in~\cite{Kolesnik2012} it follows that
\begin{equation}\label{eq:igual}
\begin{split}
\mathbb{E}_{\mu}&\left[e^{2aL^{-1}X^{\textsf{sym},v_0}(s)}\right]\\
&=e^{-\lambda s}\left(
\cosh\left(s\sqrt{\lambda^2+4a^2L^{-2}v^2_0}\right)+\frac{\lambda}{\sqrt{\lambda^2+4a^2L^{-2}v^2_0}}
\sinh\left(s\sqrt{\lambda^2+4a^2L^{-2}v^2_0}\right)
\right)
\end{split}
\end{equation}
for $4a^2L^{-2}<\lambda^2/v^2_0$, where $\cosh(\cdot)$ and $\sinh(\cdot)$ are the hyperbolic cosine and hyperbolic sine functions, respectively.
Note that both sides of 
the~\eqref{eq:igual} are analytic (holomorphic) in a strip containing the real axis, therefore the formula holds true for all real values, see corollary at the end of p.~209 in~\cite{Rudin1987}.
Since
\begin{equation}
\begin{split}
\cosh(x)=\frac{1}{2}(e^{x}+e^{-x})\leq e^{x}\quad \textrm{ and }\quad 
\sinh(x)=\frac{1}{2}(e^{x}-e^{-x})\leq e^{x}
\end{split}
\end{equation}
for all $x\geq 0$, we have 
that~\eqref{eq:igual} yields
\begin{equation}\label{eq:igual1}
\begin{split}
\mathbb{E}_{\mu}&\left[e^{2aL^{-1}X^{\textsf{sym},v_0}(s)}\right]\leq e^{-\lambda s}e^{s\sqrt{\lambda^2+4a^2L^{-2}v^2_0}}
\Big(
1+\frac{\lambda}{\sqrt{\lambda^2+4a^2L^{-2}v^2_0}}
\Big).
\end{split}
\end{equation}
The preceding inequality with the help 
of~\eqref{eq:fb} implies that
\begin{equation}\label{eq:West}
\begin{split}
\mathbb{E}_{\mu}[(W_{a,b}(L^{-1}X^{\textsf{sym},v_0}))^2]
&\leq \left(1+\frac{\lambda}{\rho}\right)
T^{-1}\int_{0}^T \rmd s \,e^{2bs}e^{-\lambda s}e^{\rho s}\\
&=\Big(\frac{\rho+\lambda}{\rho}\Big)\frac{1}{T(2b+\rho-\lambda)}(e^{T(2b+\rho-\lambda)}-1)\\
&=\Big(1+\frac{1}{\sqrt{1+4a^2L^{-2}_*}}\Big)\frac{1}{2T_*b/\lambda+\frac{4a^2T_*L^{-2}_*}{1+\sqrt{1+4a^2L^{-2}_*}}}\\
&\quad\times\left(e^{2T_*b/\lambda}\exp\left(\frac{4a^2T_*L^{-2}_*}{1+\sqrt{1+4a^2L^{-2}_*}}\right)-1\right),
\end{split}
\end{equation}
where $\rho:=\sqrt{\lambda^2+4a^2L^{-2}v^2_0}>\lambda$.
We claim that
\begin{equation}\label{eq:ex}
0<\frac{e^x-1}{x}\leq \max\{1,e^x\}\quad \textrm{ for all }\quad x\in \mathbb{R}\setminus\{0\}.
\end{equation}
Indeed, on the one hand for $x>0$ we have 
$e^x-1=\int_0^x \rmd y\, e^y $. Then by the change of variable $y\mapsto xz$ we have $e^x-1=x\int_0^1 \rmd z\, e^{xz} $. Hence, the right-hand side of~\eqref{eq:ex} holds true via monotonicity. On the other hand, the result for $x<0$ holds true by the well-known inequality
$1-e^{-(-x)}\leq -x$.

By~\eqref{eq:West} and~\eqref{eq:ex} we deduce 
\begin{equation}
\begin{split}
\mathbb{E}_{\mu}&[(W_{a,b}(L^{-1}X^{\textsf{sym},v_0}))^2]\\
&\quad\qquad\leq \left(1+\frac{1}{\sqrt{1+4a^2L^{-2}_*}}\right)\max\left\{1,\exp\left({2T_*b/\lambda+\frac{4a^2T_*L^{-2}_*}{1+\sqrt{1+4a^2L^{-2}_*}}}\right)\right\}.
\end{split}
\end{equation}
The proof is complete.
\qedw

\subsection{\textbf{Proof of Corollary~\ref{corB}}}\label{sub:pcorB}\hfill

\noindent
The proof of Corollary~\ref{corB} follows step by step from the proof of Corollary~\ref{corA}.
\qedw

\section{\textbf{Numerical experiments}}\label{sec:simulations}
To assess the sharpness of the bounds obtained in Section~\ref{sub:results}, we perform a small numerical experiment to numerically evaluate the error
\begin{equation}\label{eq:num_error}
\mathbb{E}_{\nu}[F_{a,b}(B)] - \mathbb{E}_{\mu}[F_{a,b}(L^{-1}X^{\textsf{sym},v_0})],
\end{equation}
where
\begin{equation}
F_{a,b}(X):=g\left(\frac{1}{T}\int_{0}^T \rmd s \,f\big(e^{aX(s)+bs}\big)\right)
\end{equation}
with
\begin{align*}
L &:=  1 =: T &  f(x) & := x, & a & := \sigma \sqrt{\lambda}, & b & := -\frac{1}{2}\sigma^2, & v_0 & := 1, & g(x) := \max\{x - K, 0 \}
\end{align*}
fixed and
\begin{align*}
K & \in \{0.7, 1, 1.3\}, &
\sigma & \in \{0.3,0.5,0.7\}, &
\lambda  \in \{2.5,5,7.5,\ldots,100\}
\end{align*}
variable. The choice of $a$ makes the parametrization consistent with the Kac regime as $\lambda \to \infty$.
The setting models Monte Carlo pricing of an arithmetic Asian call option in finance, but our interest here is primarily computational, and we do not put it forward as a competitive method to price the option.

\begin{figure}[t!]
\begin{center}
\includegraphics[clip=TRUE, trim=0.25cm 1.4cm 0.8cm 1cm, page=2, scale=0.484]{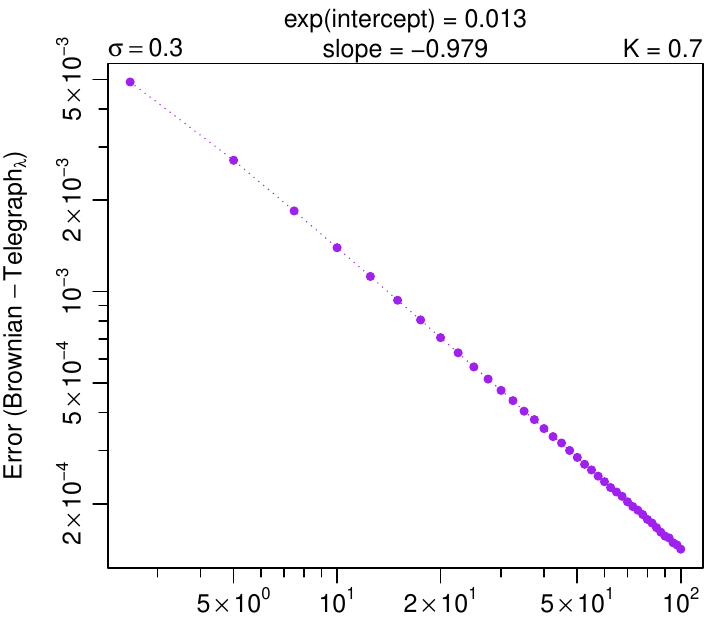}
\includegraphics[clip=TRUE, trim=1.2cm 1.4cm 0.8cm 1cm, page=4, scale=0.484]{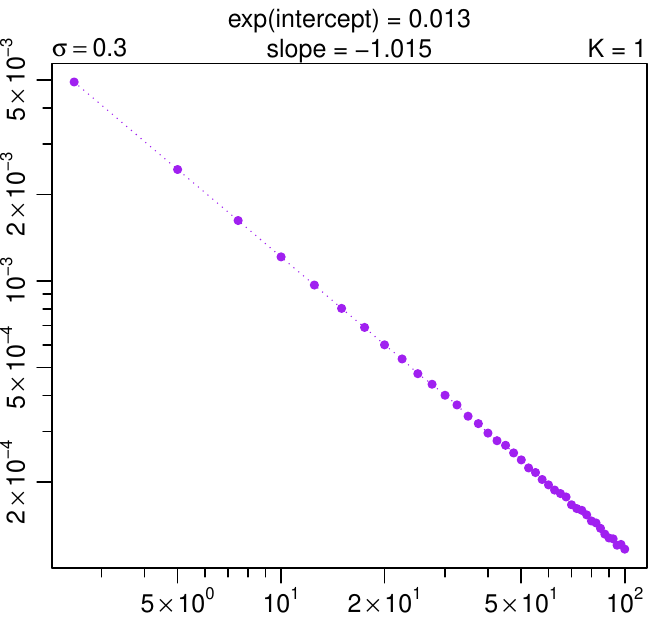}
\includegraphics[clip=TRUE, trim=1.2cm 1.4cm 0.8cm 1cm, page=6, scale=0.484]{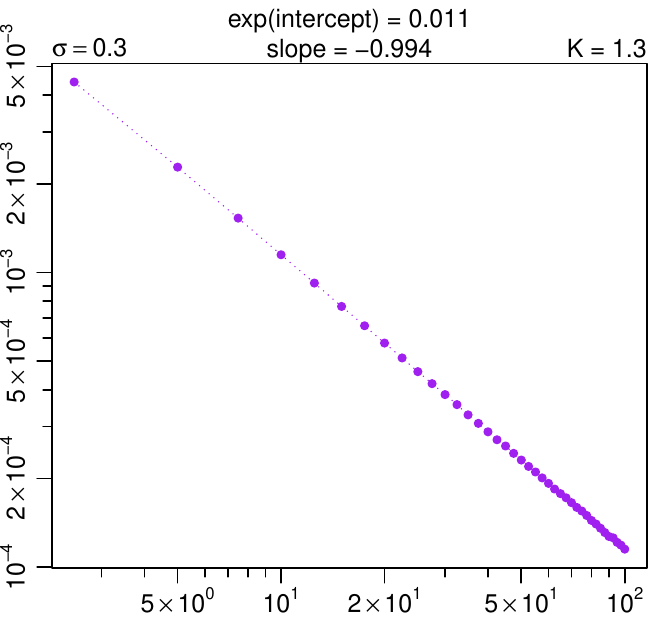} \\
\includegraphics[clip=TRUE, trim=0.25cm 1.4cm 0.8cm 1cm, page=8, scale=0.484]{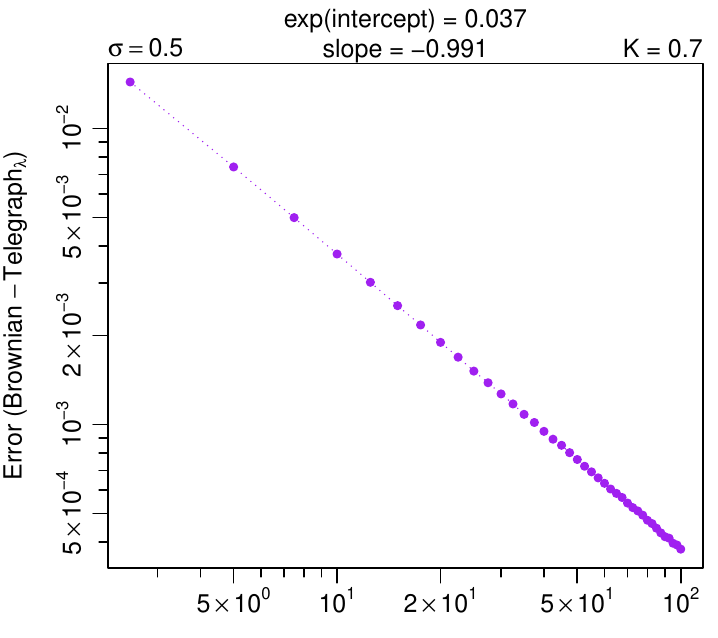}
\includegraphics[clip=TRUE, trim=1.2cm 1.4cm 0.8cm 1cm, page=10, scale=0.484]{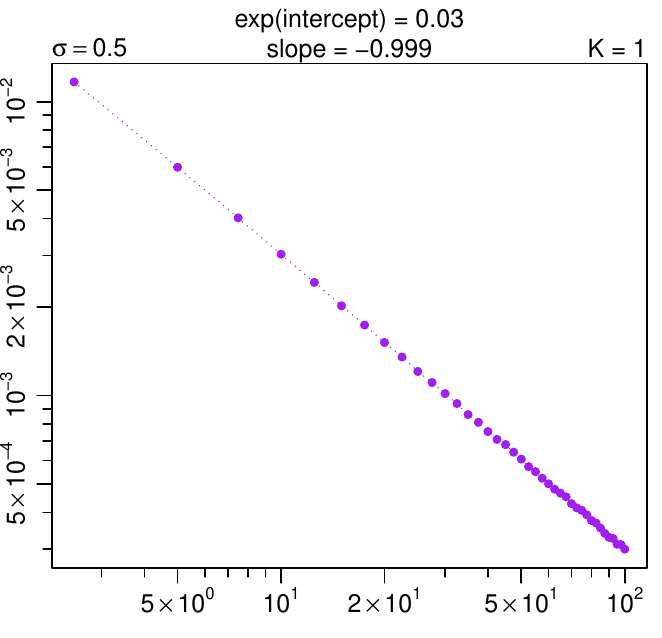}
\includegraphics[clip=TRUE, trim=1.2cm 1.4cm 0.8cm 1cm, page=12, scale=0.484]{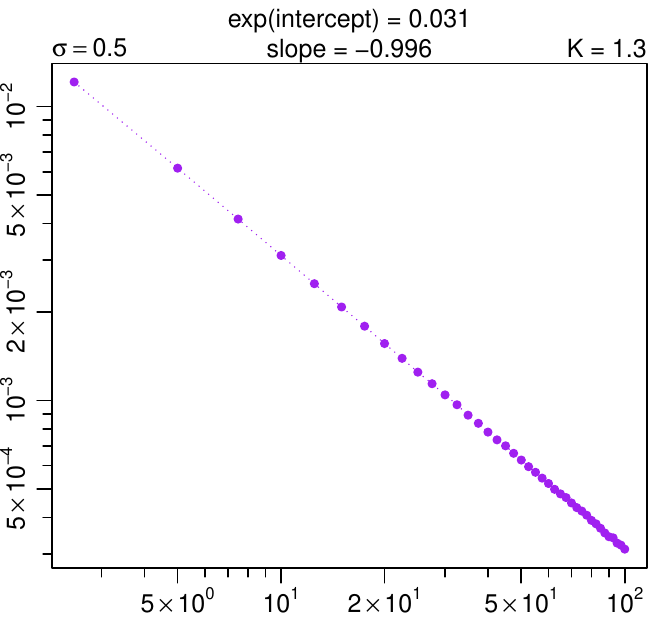} \\
\includegraphics[clip=TRUE, trim=0.25cm 0.8cm 0.8cm 1cm, page=14, scale=0.484]{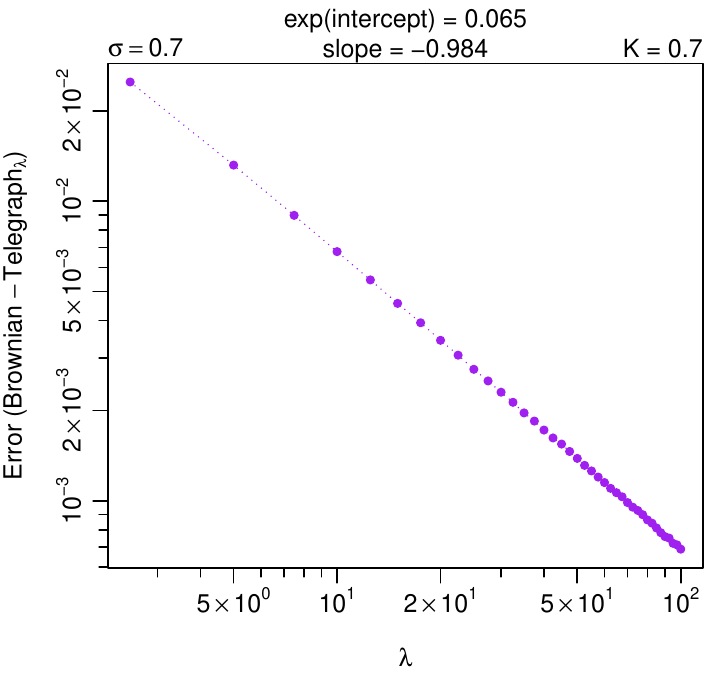}
\includegraphics[clip=TRUE, trim=1.2cm 0.8cm 0.8cm 1cm, page=16, scale=0.484]{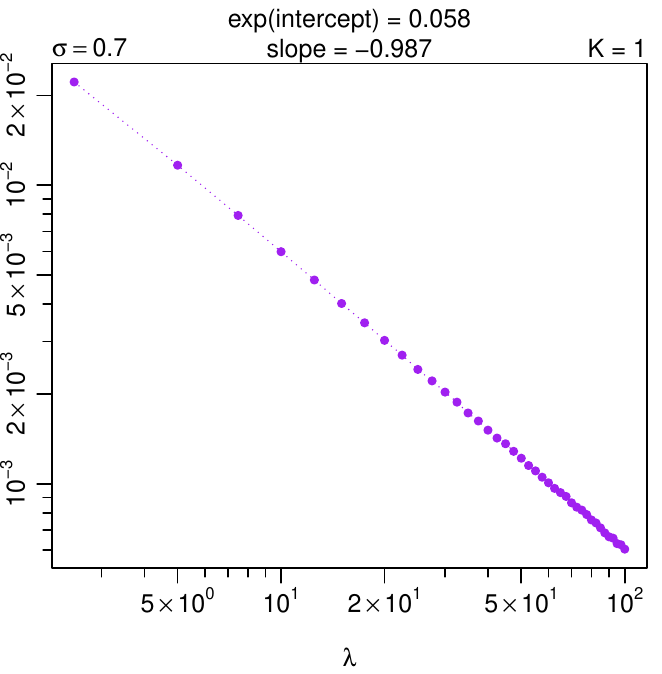}
\includegraphics[clip=TRUE, trim=1.2cm 0.8cm 0.8cm 1cm, page=18, scale=0.484]{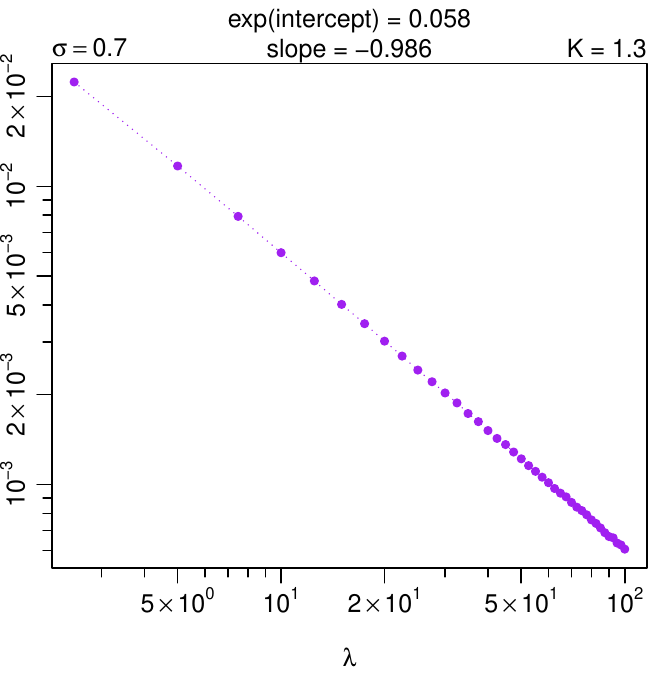}
\end{center}
\caption{Numerical estimates of the approximation 
error~\eqref{eq:num_error} as a function of $\lambda$ and the intercept and slope estimates of an ordinary least squares regression of the natural logarithm of~\eqref{eq:num_error} on $\ln \lambda$.}\label{fig:simu_res}
\end{figure}

To exactly simulate the random variable $F_{a,b}(X^{\textsf{sym},1})$ in closed form in the case $f(x) = x$, we proceed as follows. Note first that $N(T) + N^*(T) \sim \mathrm{Poisson}(2\lambda)$ and that, conditional on $N(T) + N^*(T) = n$ for some $n \in \mathbb{N}$, the collated and ordered jump times $T_1 \leq \cdots \leq T_n$ of $N$ and $N^*$ are i.i.d.\ draws from the $\mathrm{Uniform}(0,T)$ distribution sorted in increasing order. (Note that the probability of ties is zero.) Subsequently, we compute
\begin{equation*}
F_{a,b}(X^{\textsf{sym},1}) = g\left(\frac{1}{T}\int_0^T \rmd s \,e^{aX(s)+bs}\right) = g\left(\sum_{i=1}^{n+1} \int_{T_{i-1}}^{T_i} \rmd s \,e^{aX(s)+bs}\right) = g\left(\sum_{i=1}^{n+1} A_i\right)\,,
\end{equation*}
where $T_0 := 0$, $T_{n+1}:= T$ and
\begin{equation*}
A_i := \frac{1}{(-1)^{i-1}a + b} \big(e^{ ((-1)^{i-1}a+b) (T_i - T_{i-1})}-1\big) e^{a \sum_{j=1}^{i-1} (-1)^{j-1} (T_j-T_{j-1})+b T_{i-1}},
\end{equation*}
interpreting the empty sum in the case $i=1$ as zero.

To estimate $\mathbb{E}_{\mu}[F_{a,b}(X^{\textsf{sym},1})]$, we average over $10^9$ realizations of $F_{a,b}(X^{\textsf{sym},1})$ for each combination of the values of the variable parameters $K$, $\sigma$ and $\lambda$. Moreover, we simulate the random variable $F_{a,b}(B)$ by sampling the standard Brownian motion $B$ at $10^4$ equidistant time points covering $[0,1]$ and discretizing the integral in $F_{a,b}(B)$ using Riemann sums. Finally, we estimate $\mathbb{E}_{\nu}[F_{a,b}(B)]$ by averaging over $10^9$ such realizations. The C++ implementations of the simulation routines, to be run in conjunction with R~\cite{R2025} and 
Rcpp~\cite{Eddelbuettel2025}, are provided in \url{https://github.com/mspakkanen/telegraph}.

The results of this experiment are shown in Figure~\ref{fig:simu_res}. To assess the rate of convergence of the error~\eqref{eq:num_error} to zero as $\lambda \rightarrow \infty$, we regress its natural logarithm on $\ln (\lambda)$ and report the intercept and slope estimates. The results overwhelmingly suggest that the error behaves like $\lambda^{\alpha}$ with exponent $\alpha \approx -1$ as $\lambda \rightarrow \infty$, 
the behavior being uniform over the different values of $\sigma$ and $K$.
This asymptotic behavior is also consistent with our earlier theoretical results where $\alpha \approx -1/4 > -1$ for $\lambda \to \infty$,
see Remark~2.8 in~\cite{BarreraLukkarinen2023}.
The  bounds in~\cite{BarreraLukkarinen2023} are uniform in $\lambda$ and indeed their proof suggests that the main errors come from the regime of small values of the rate $\lambda$.

\section*{\textbf{Statements and declarations}}

\noindent
\textbf{Acknowledgments.} Gerardo Barrera would like to express his gratitude to University of Helsinki (Finland) and Instituto Superior T\'ecnico, (Portugal) for all the facilities used along the realization of this work. 
The authors are indebted with professor Dario Gasbarra (University of Vaasa, Finland) for valuable discussions  during the preparation of this manuscript.
They also acknowledge computational resources and support provided by the Imperial College Research Computing Service (\url{http://doi.org/10.14469/hpc/2232}).

\noindent
\textbf{Funding.}
The work has been supported by the Academy of Finland, via an Academy project (project No.\ 339228) and the \textit{Finnish Centre of Excellence in Randomness and Structures} (projects No.\  346306, 364213).
The research of Gerardo Barrera is partially supported by Horizon Europe Marie Sk\l{}odowska-Curie Actions Staff Exchanges (project no.\ 101183168 -- LiBERA).

\noindent
\textbf{Availability of data and material.} 
The results generated in the numerical experiment are available upon request.

\noindent
\textbf{Conflict of interests.} The authors declare that they have no conflict of interest.

\noindent
\textbf{Ethical approval.} Not applicable.

\noindent
\textbf{Authors' contributions.}
All authors have contributed equally to the paper.

\noindent
\textbf{Disclaimer.} 
Funded by the European Union. Views and opinions expressed are however those of the author(s) only and do not necessarily reflect those of the European Union or the European Education and Culture Executive Agency (EACEA). Neither the European Union nor EACEA can be held responsible for them.

\bibliographystyle{amsplain}

\end{document}